\documentclass[11pt]{article}

\usepackage{amsmath, amsfonts, amssymb,latexsym,wasysym,graphicx}
\input epsf.sty

\newtheorem{Theorem}{Theorem}[section]
\newtheorem{Lemma}{Lemma}[section]

\newtheorem{Definition}[Lemma]{Definition}

\newcommand{\BEQ}{\begin{equation}}     
\newcommand{\BEA}{\begin{eqnarray}}
\newcommand{\BD}{\begin{displaymath}}
\newcommand{\EEQ}{\end{equation}}       
\newcommand{\EEA}{\end{eqnarray}}
\newcommand{\ED}{\end{displaymath}}

\newcommand{\del}{\delta}

\newcommand{\eps}{\varepsilon}          

\newcommand{\R}{\mathbb{R}}

\newcommand{\Sk}{{\mathrm{Sk}}}

\def\esper{{\mathbb{E}}}
\def\T{{\mathbb{T}}}
\def\Var{{\mathrm{Var}}}

\newcommand{\eop}{\hfill $\Box$}        

\newcommand{\II}{{\rm i}}               
\renewcommand{\Re}{{\rm Re\ }}          
\newcommand{\half}{{1\over 2}}          
\renewcommand{\vec}[1]{\boldsymbol{#1}} 

\catcode`\@=11
\def\numberbysection{\@addtoreset{equation}{section}
        \def\theequation{\thesection.\arabic{equation}}}
\numberbysection

\begin{document}

\vspace*{1.5cm}
\begin{center}
{\Large \bf A L\'evy area by Fourier normal ordering  for
multidimensional fractional Brownian motion with small Hurst index
}
\end{center}

\vspace{2mm}
\begin{center}
{\bf  J\'er\'emie Unterberger}
\end{center}

\vspace{2mm}
\begin{quote}

\renewcommand{\baselinestretch}{1.0}
\footnotesize
{The main tool for stochastic calculus with respect to a multidimensional process $B$ with 
small H\"older regularity index is rough path theory. Once $B$ has been lifted to a rough
path, a stochastic calculus -- as well as solutions to stochastic differential equations driven
by $B$ -- follow by standard arguments. 

Although such a lift has been proved to exist by abstract arguments \cite{LyoVic07}, a first general, explicit construction
has been proposed in \cite{Unt09,Unt09bis} under the name of Fourier normal ordering.

The purpose of this short note is to convey the main ideas of the Fourier normal ordering
method in the particular case of the iterated integrals of lowest order of fractional Brownian
motion with arbitrary Hurst index.
 }
\end{quote}

\vspace{4mm}


\newpage


\section{Introduction}


The (two-sided) fractional Brownian motion $t\to B_t$, $t\in\R$  (fBm for short) with Hurst exponent $\alpha$, $\alpha\in(0,1)$, defined as the centered Gaussian process with covariance
\BEQ \esper[B_s B_t]=\half (|s|^{2\alpha}+|t|^{2\alpha}-|t-s|^{2\alpha}), \EEQ
is a natural generalization in the class of Gaussian processes of
the usual Brownian motion (which is the case $\alpha=\half$), in the sense that it exhibits two fundamental properties shared with Brownian motion, namely,
it has stationary increments, viz. $\esper[(B_t-B_s)(B_u-B_v)]=\esper[(B_{t+a}-B_{s+a})(B_{u+a}-B_{v+a})]$ for
every $a,s,t,u,v\in\R$, and it is self-similar, viz. 
\BEQ 
\forall \lambda>0, \quad (B_{\lambda t}, t\in\R) \overset{(law)}{=} (\lambda^{\alpha} B_t,
t\in\R).
\EEQ 
One may also define  a $d$-dimensional 
vector Gaussian process (called: {\it $d$-dimensional fractional Brownian motion}) by setting $B_t=(B_t(1),\ldots,B_t(d))$ where $(B_t(i),t\in\R)_{i=1,\ldots,d}$ are $d$ independent (scalar) fractional Brownian motions.

Its theoretical interest lies in particular in the fact that it is (up to normalization) the only Gaussian process satisfying these two properties.

A standard application of Kolmogorov's theorem shows that fBm has a version with
$\alpha^-$-H\"older continuous (i.e.  $\kappa$-H\"older continuous for every $\kappa<\alpha$) paths.
In particular, fBm with small Hurst parameter $\alpha$ is a natural, simple  model for continuous but
very irregular processes.

There has been a widespread interest during the past ten years in constructing a stochastic integration theory
with respect to fBm and solving stochastic differential equations driven by fBm, see for instance
\cite{LLQ02,GraNouRusVal04,CheNua05,RusVal93,RusVal00}. The multi-dimensional case is
very different from the one-dimensional case. When one tries to integrate for instance a stochastic differential
equation driven by a two-dimensional fBm $B=(B(1),B(2))$ by using any kind of Picard iteration scheme, one
encounters very soon the problem of defining the L\'evy area of $B$ which is the antisymmetric part
of ${\cal A}_{ts}:=\int_s^t dB_{t_1}(1) \int_s^{t_1} dB_{t_2}(2)$. This is the simplest occurrence
of iterated integrals $\vec{B}^k_{ts}(i_1,\ldots,i_k):=\int_s^t dB_{t_1}(i_1)\ldots \int_s^{t_{k-1}} dB_{t_k}(i_k)$,
$i_1,\ldots,i_k\le d$
for $d$-dimensional fBm $B=(B(1),\ldots,B(d))$ 
which lie at the heart of the rough path theory due to T. Lyons, see \cite{Lyo98,LyoQia02}. 
An alternative construction has been given in \cite{Gu} under the name of 'algebraic
rough path theory', which we now propose to 
describe  briefly.

 Assume
$\Gamma_t=(\Gamma_t(1),\ldots,\Gamma_t(d))$ is some non-smooth $d$-dimensional path
which is $\alpha$-H\"older continuous.  Integrals such as
$\int f_1(\Gamma_t)d\Gamma_t(1)+\ldots+f_d(\Gamma_t)d\Gamma_t(d)$ do not make sense a priori because
$\Gamma$ is not differentiable (Young's integral \cite{Lej} works for $\alpha>\half$ but not beyond). In order to define
the integration of a differential form along $\Gamma$, it is enough to define a {\it truncated multiplicative functional}
or {\it geometric rough path} $(\vec{\Gamma}^{1},\ldots,\vec{\Gamma}^{\lfloor 1/\alpha \rfloor})$ lying above $\Gamma$,  
 $\lfloor  1/\alpha \rfloor$=entire part of $1/\alpha$,
  where $\vec{\Gamma}^{1}_{ts}=(\del\Gamma)_{ts}:=\Gamma_t-\Gamma_s$ is the {\it increment}
of $\Gamma$ between $s$ and $t$, and
each
$\vec{\Gamma}^k=(\vec{\Gamma}^k(i_1,\ldots,i_k))_{1\le i_1,\ldots,i_k\le d}$, $k\ge 2$ 
is a {\it substitute} for the iterated integrals $\int_s^t d\Gamma_{t_1}(i_1)\int_s^{t_1} d\Gamma_{t_2}(i_2)
\ldots \int_{s}^{t_{k-1}} d\Gamma_{t_k}(i_k)$ with the following three properties:

\begin{itemize}
\item[(i)] (\it H\"older continuity)
 each component of $\vec{\Gamma}^k$ is $k\alpha^-$-H\"older continuous, that is
to say, $k\kappa$-H\"older for every $\kappa<\alpha$;

\item[(ii)] ({\it multiplicative or Chen property}) letting $\del{\bf\Gamma}^k_{tus}:=
{\bf\Gamma}_{ts}^k-{\bf\Gamma}^k_{tu}-{\bf \Gamma}^k_{us}$, one requires
\BEQ
 \del\vec{\Gamma}^k_{tus}(i_1,\ldots,i_k) = \sum_{k_1+k_2=k} \vec{\Gamma}_{tu}^{k_1}(i_1,\ldots,i_{k_1}) \vec{\Gamma}_{us}^{k_2}(i_{k_1+1},\ldots,i_k). \label{eq:0:x}
\EEQ

\item[(iii)] ({\it geometric or shuffle property}) 
\BEQ  {\bf \Gamma}^{n_1}_{ts}(i_1,\ldots,i_{n_1}) {\bf \Gamma}^{n_2}_{ts}(j_1,\ldots,j_{n_2}) 
 = 
\sum_{\vec{k}\in {\mathrm{Sh}}(\vec{i},\vec{j})} {\bf \Gamma}^{n_1+n_2}_{ts}(k_1,\ldots,k_{n_1+n_2})  \label{eq:0:geo} \EEQ
where ${\mathrm{Sh}}(\vec{i},\vec{j})$ is the subset of permutations of $i_1,\ldots,i_{n_1},j_1,\ldots,j_{n_2}$
which do not change the orderings of $(i_1,\ldots,i_{n_1})$ and $(j_1,\ldots,j_{n_2})$.
\end{itemize}

 Then there
is a standard procedure which allows to define out of these data iterated integrals of any order and
to solve differential equations driven by $\Gamma$.

The multiplicativity property (\ref{eq:0:x}) and the geometric property (\ref{eq:0:geo}) are satisfied by smooth paths, as can be checked by direct computation. So the most
natural way to construct such a multiplicative functional is to start from some smooth approximation
$\Gamma^{\eps}$, $\eps\overset{>}{\to} 0$ of $\Gamma$ such that each iterated
 integral $\vec{\Gamma}^{k,\eps}_{ts}(i_1,\ldots,i_k)$, $k\le \lfloor 1/\alpha \rfloor$
 converges in the 
$k\kappa$-H\"older norm for every $\kappa<\alpha$.

This general scheme has been applied to fBm in a paper by L. Coutin and Z. Qian \cite{CQ02} and later in a paper
by the author \cite{Unt08}, using different schemes of approximation of $B$ by a family of Gaussian processes $B^{\eps}$ (living on the same probability
space) with
$\eps\to 0$.
 In both cases, the variance of the L\'evy area has been proved
to diverge in the limit $\eps\to 0$ when $\alpha\le 1/4$.

\medskip

Let us explain briefly our construction for the second-order iterated
integral ${\bf B}^2_{ts}(i_1,i_2)$ (by abuse of language, we shall
call this object a L\'evy area, although the L\'evy area is usually
defined as  the
corresponding antisymmetrized quantity).

\medskip

Let $\alpha\in(0,\half)$ (or even  $\alpha\in(0,1/4)$), and consider
the natural iterated integral ${\bf B}^{2,\eps}_{ts}(i_1,i_2)$ for some family of
approximations $B^{\eps}$ of $B$. Assume ${\bf Z}^{\eps}(i_1,i_2)$
 is some a.s. $2\alpha^-$-H\"older random function
 (living on the same probability
space as $B^{\eps}$), antisymmetric in $(i_1,i_2)$. Then 
\BEQ {\cal R} {\bf B}^{2,\eps}_{ts}(i_1,i_2):={\bf B}^{2,\eps}_{ts}(i_1,i_2)+
\del {\bf Z}^{\eps}_{ts}(i_1,i_2)=
{\bf B}^{2,\eps}_{ts}(i_1,i_2)+ \left({\bf Z}^{\eps}_t(i_1,i_2)-
{\bf Z}^{\eps}_s(i_1,i_2) \right) \EEQ
satisfies properties (i), (ii) and (iii). The multiplicativity property
is preserved because $\del$ acting on an increment vanishes \cite{Gu}, 
\BEA \del(\del {\bf Z}^{\eps}(i_1,i_2))_{tus} &=&  \left ({\bf Z}^{\eps}_t(i_1,i_2)-{\bf Z}^{\eps}_s(i_1,i_2)\right)-
\left ({\bf Z}^{\eps}_t(i_1,i_2)-{\bf Z}^{\eps}_u(i_1,i_2)\right) \nonumber\\
&-&
\left ({\bf Z}^{\eps}_u(i_1,i_2)-{\bf Z}^{\eps}_s(i_1,i_2)\right) =0, \label{eq:0:deldel} 
\EEA
and the geometric property is also preserved because ${\bf Z}^{\eps}(i_1,i_2)+{\bf Z}^{\eps}(i_2,i_1)=0$.

The function ${\bf Z}^{\eps}$
 may be seen as a {\it counterterm}. Now (see section 1) ${\bf Z}^{\eps}$
may be chosen so as to make the {\it regularized L\'evy area} ${\cal R}
{\bf B}^{2,\eps}_{ts}(i_1,i_2)$ converge for every $\kappa<\alpha$ and $T>0$
 in $L^2(\Omega,{\cal C}^{2\kappa}([-T,T]))$
to a finite L\'evy area for $B$, where ${\cal C}^{2\kappa}([-T,T])$ is
the Banach space of $2\kappa$-H\"older 1-increments on $[-T,T]$ in the sense of Gubinelli \cite{Gu},
equipped with the H\"older norm $||f||_{2\kappa,T}=\sup_{s,t\in[-T,T]} \frac{|f_{t,s}|}{|t-s|^{2\kappa}}$.
  More
precisely, ${\bf Z}^{\eps}$ may be chosen in the second chaos of $B$. One
may prove bounds of the type 
\BEQ \esper  {\cal R}{\bf B}^{2,\eps}_{ts}(i_1,i_2)^2
\le C|t-s|^{4\alpha}, \quad \esper \left(  {\cal R}{\bf B}^{2,\eps}_{ts}(i_1,i_2)-{\cal R}
{\bf B}^{2,\eta}_{ts}(i_1,i_2) \right)^2 \le C|\eps-\eta|^{2\alpha} \EEQ
from which (see \cite{Unt09}, Proposition 1.5) the convergence in 
   $L^2(\Omega,{\cal C}^{2\kappa}([-T,T]))$ may be proved by using standard arguments, in 
particular the Garsia-Rodemich-Rumsey lemma \cite{Ga}.

\medskip

The first section concerns the construction of the counterterm ${\bf Z}^{\eps}$ and
the regularization of the L\'evy area. The main tool is  Fourier transform. Splitting iterated
integrals into {\it increment/boundary terms} (see below) and reordering Fourier
components by
Fubini's theorem in such a way that {\it innermost integrals bear highest Fourier frequencies}
ensures the proper H\"older regularity  for every {\it boundary term}, hence the name of {\it Fourier normal
ordering}. Then {\it increment terms}
are regularized by introducing an adequate cut in the Fourier domain of integration. The
multiplicative rule is preserved by doing this (see eq. (\ref{eq:0:deldel})).

\medskip

 The sketchy generalization to higher-order iterated integrals in section 2 involves
 tree combinatorics linked with
the {\it Hopf algebra} structure of the algebra of decorated rooted trees, as defined by A. Connes
and D. Kreimer in a series of papers 
\cite{ConKre98,ConKre00,ConKre01} devoted to a mathematical axiomatization of the BPHZ (Bogolioubov and coauthors)
algorithm of renormalization of Feynman graphs in quantum field theory \cite{Hepp}. Our algorithm has nothing to do
with BPHZ though, since it is based on a somewhat arbitrary but convenient 
regularization of {\it skeleton integrals} (higher-order
generalizations of the increment term of the L\'evy area), which are {\it 
tree iterated integrals} of a particular type. 
 The proof of the multiplicative/Chen and geometric/shuffle property in \cite{Unt09bis} relies on Hopf algebra computations though.
A BPHZ or dimensional renormalization scheme could be applied to skeleton integrals instead of a blunt regularization (work
in progress). This algebraic approach has proved to be useful in a variety of contexts ranging from
numerical analysis to quantum chromodynamics or the study of polylogarithms or multi-zeta functions, see
for instance \cite{Kre99,Mur2,Wal00}.


\section{Definition of a L\'evy area}


Recall that fractional Brownian motion $B$ may be defined via the harmonizable representation \cite{SamoTaq}
\BEQ B_t=c_{\alpha} \int_{\R}  |\xi|^{\half-\alpha}  \frac{e^{\II t\xi}-1}{\II \xi}
\ W(d\xi) \EEQ
where $(W_{\xi},\xi\ge 0)$ is a complex Brownian motion extended to $\R$
 by setting $W_{-\xi}=-\overline{W}_{\xi}$ $(\xi\ge 0)$, and 
$c_{\alpha}=\half \sqrt{-\frac{\alpha}{\cos\pi\alpha\Gamma(-2\alpha)}}$.

We shall use the following approximation of $B$ by a family of centered Gaussian
processes $(B^{\eps},\eps>0)$ living in the first chaos of $B$.

\begin{Definition}[approximation $B^{\eps}$]

Let, for $\eps>0$,

\BEQ B_t^{\eps}=c_{\alpha} \int_{\R} e^{-\eps|\xi|} |\xi|^{\half-\alpha} 
 \frac{e^{\II t\xi}-1}{\II \xi}\ W(d\xi). \EEQ

\end{Definition}

The process $B^{\eps}$ is easily seen to have a.s. smooth paths. The infinitesimal covariance
$\esper (B^{\eps})'_s (B^{\eps})'_t$ may be computed explicitly using the Fourier transform
\cite{Erd54}
\BEQ  {\cal F}K^{',-}_{\eps}(\xi)=\frac{1}{\sqrt{2\pi}} \int_{\R} K^{',-}_{\eps}(x) e^{-\II x\xi} dx
=-\frac{\pi\alpha}{2\cos\pi\alpha\Gamma(-2\alpha)} e^{-2\eps |\xi|} |\xi|^{1-2\alpha} 
{\bf 1}_{|\xi|>0},\EEQ 
where $K^{',-}_{\eps}(s-t):=\frac{\alpha(1-2\alpha)}{2\cos\pi\alpha} (-\II(s-t)+2\eta)^{2\alpha-2}$.
By taking the real part of these expressions, one finds that $B^{\eps}$ has the same law as
the analytic approximation of $B$ defined in \cite{Unt08}, namely, $B^{\eps}=
\Gamma_{t+\II\eps}+\Gamma_{t-\II\eps}=2\Re\Gamma_{t+\II\eps}$, where $\Gamma$ is the
analytic fractional Brownian motion (see also \cite{TinUnt}).

\medskip

This second-order case is too simple to capture the combinatorial features of Fourier normal
ordering. On the other hand, the proof of convergence for the regularized L\'evy area (after
substraction of the counterterm) is short, and the proof for tree integrals in the
general case may be considered as a generalization.

Fourier normal ordering is the combination of (i) a number of equivalent splittings
 of the iterated integrals related by Fubini's theorem; (ii) a splitting of the Fourier
domain of integration into a number of disjoint Fourier domains; (iii)
an appropriate choice of splitting of the iterated integrals on each  Fourier domain.

 Let us first write down the two equivalent splittings in the case of the L\'evy area. The components
$B(1),B(2)$ are assumed to be constructed from i.i.d. Brownian motions $W(1)$, $W(2)$ via the above
harmonizable representation. Note that ${\bf B}^{2,\eps}_{ts}(j,j)=\half(B_t^{\eps}(j)-B_s^{\eps}(j))^2$,
$j=1,2$
requires no regularization, hence the only problematic second-order iterated integrals are the
mixed integrals  ${\bf B}^{2,\eps}(i,j)$, $i\not=j$, say, $i=1,j=2$.

\begin{Lemma}

Let ${\bf B}^{2,\eps}_{ts}(1,2)=\int_s^t dB_{u_1}^{\eps}(1)\int_s^{u_1} dB^{\eps}_{u_2}(2).$ Then:

\BEQ {\bf B}^{2,\eps}_{ts}(1,2)=c_{\alpha}^2\int_{-\infty}^{+\infty} dW_{\xi_1}(1) \int_{-\infty}^{+\infty}
dW_{\xi_2}(2) e^{-\eps(|\xi_1|+|\xi_2|)} |\xi_1 \xi_2|^{\half-\alpha} I_{ts}(\xi_1,\xi_2) \EEQ
where
\BEQ I_{ts}(\xi_1,\xi_2):=
\int_s^t e^{\II u_1\xi_1} du_1 \int_s^{u_1} e^{\II u_2\xi_2} du_2.\EEQ

\end{Lemma}

The proof is straightforward.

There are two apparently equivalent ways of splitting the integral $I_{ts}(\xi_1,\xi_2)$ into
an {\bf increment term}, $(\del J)_{ts}=J_t-J_s$, and a {\bf boundary term} denoted by the symbol $\partial$:

\begin{itemize}
\item[(i)] either writing $I_{ts}(\xi_1,\xi_2)$ as $(G_t^+-G_s^+)(\xi_1,\xi_2)+I_{ts}^+(\xi_1,\xi_2)(\partial)$, where (provided $\xi_1+\xi_2\not=0$)
\BEQ G_u^+(\xi_1,\xi_2)=\frac{e^{\II u(\xi_1+\xi_2)}}{[\II(\xi_1+\xi_2)][\II\xi_2]}, \qquad
I_{ts}^+(\xi_1,\xi_2)(\partial)=-\frac{e^{\II s\xi_2}}{\II \xi_2} \ .\ \frac{e^{\II t\xi_1}
-e^{\II s\xi_1}}{\II \xi_1};  \label{eq:G+I+} \EEQ
\item[(ii)] or (using Fubini's theorem) 
\BEA I_{ts}(\xi_1,\xi_2) &=& \int_s^t e^{\II u_2\xi_2} du_2 \int_{u_2}^t e^{\II u_1\xi_1} du_1 \nonumber\\
&=& (G_t^- -G_s^-)(\xi_2,\xi_1)+
I_{ts}^-(\xi_2,\xi_1)(\partial), \EEA
 where
\BEQ G_u^-(\xi_2,\xi_1)=-\frac{e^{\II u(\xi_1+\xi_2)}}{[\II(\xi_1+\xi_2)][\II\xi_1]}, \qquad
I_{ts}^-(\xi_2,\xi_1)(\partial)=\frac{e^{\II t\xi_1}}{\II \xi_1} \ .\ \frac{e^{\II t\xi_2}
-e^{\II s\xi_2}}{\II \xi_2}. \EEQ
\end{itemize}

In either case, the inner integral $\int_x^u e^{\II u'\xi} du'$, $(u,x,\xi)=(u_1,s,\xi_2)$ or $(u_2,t,\xi_1)$ has been formally
decomposed as $\int^u e^{\II u'\xi} du'-\int^x e^{\II u'\xi} du'=\frac{e^{\II u\xi}}{\II \xi}
-\frac{e^{\II x\xi}}{\II\xi}$, which introduces an apparent infra-red divergence, since the 
single terms $\frac{e^{\II u\xi}}{\II \xi}$, $\frac{e^{\II x\xi}}{\II \xi}$ diverge when $\xi\to 0$.
\footnote{Formally $\int^u e^{\II u\xi} du=\frac{e^{\II u\xi}}{\II \xi}=\int_{\pm\II\infty}^u
e^{\II u\xi} du$ depending on the sign of $\xi$.} Boundary terms come from the contribution of
$\frac{e^{\II x\xi}}{\II \xi}$, $x=s$ or $t$, and are {\it not} increments.

\medskip

The idea of Fourier normal ordering is that {\bf innermost integrals should bear highest Fourier frequencies} in order to get correct H\"older estimates {\it separately} for the increment and the boundary
term. Namely, using for instance the decomposition eq. (\ref{eq:G+I+}) for arbitrary values of $\xi_1,\xi_2$ yields
a boundary term $-B_s^{\eps}(2)(B_t^{\eps}(1)-B_s^{\eps}(1))$ which is obviously only $\alpha^-$-H\"older, and not $2\alpha^-$-H\"older.

We shall say that a function of two arguments
 $(t,s)\to {\bf X}^{\eps}_{ts}$ ($\eps>0$) in the second chaos of $B$ is {\it uniformly $2\alpha^-$-H\"older
in $\eps$} if
\begin{itemize}
\item[(i)] $\esper ({\bf X}^{\eps}_{ts})^2 \le C|t-s|^{4\alpha}$ for a constant $C$ which is
independent of $\eps$, and furthermore  
\item[(ii)] one has the following {\it rate of convergence}
  $\esper \left( {\bf X}^{\eps}_{ts}-{\bf X}^{\eta}_{ts}\right)^2 
\le C |\eps-\eta|^{2\alpha}$ when $\eps,\eta\to 0$.
\end{itemize}

 This implies (by the arguments given
in the Introduction) that ${\bf X}^{\eps}_{ts}$
converges in $L^2(\Omega;{\cal C}^{2\kappa}([-T,T]))$ for every $T>0$ and $\kappa<\alpha$.

\begin{Lemma} \label{lem:B2}

\begin{itemize}
\item[(i)] ({\bf boundary term})  For every $\alpha\in(0,\half)$,
\BEQ {\bf B}_{ts}^{2,\eps,+}(1,2)(\partial):=c_{\alpha}^2\int\int_{|\xi_1|\le |\xi_2|} dW_{\xi_1}(1)
dW_{\xi_2}(2) e^{-\eps(|\xi_1|+|\xi_2|)} |\xi_1 \xi_2|^{\half-\alpha} I_{ts}^+(\xi_1,\xi_2)(\partial)
\EEQ
and
\BEQ {\bf B}_{ts}^{2,\eps,-}(1,2)(\partial):=c_{\alpha}^2\int\int_{|\xi_2|\le |\xi_1|} dW_{\xi_1}(1)
dW_{\xi_2}(2) e^{-\eps(|\xi_1|+|\xi_2|)} |\xi_1 \xi_2|^{\half-\alpha} I_{ts}^-(\xi_2,\xi_1)(\partial)
\EEQ
are $2\alpha^-$-H\"older  uniformly in $\eps$.
\item[(ii)] ({\bf increment term}) For every $\alpha\in(1/4,1/2)$, the functions
\BEQ  \del G_{ts}^{2,\eps,+}(1,2):=c_{\alpha}^2\int\int_{|\xi_1|\le |\xi_2|}  e^{-\eps(|\xi_1|+|\xi_2|)} |\xi_1 \xi_2|^{\half-\alpha} \left(G_t^+(\xi_1,\xi_2)-G_s^+(\xi_1,\xi_2)\right)\EEQ
and
\BEQ \del G_{ts}^{2,\eps,-}(1,2):=c_{\alpha}^2\int\int_{|\xi_2|\le |\xi_1|}  e^{-\eps(|\xi_1|+|\xi_2|)} |\xi_1 \xi_2|^{\half-\alpha} \left(G_t^-(\xi_2,\xi_1)-G_s^-(\xi_2,\xi_1)\right) \EEQ
satisfy the first estimates (i) $\esper |\del G^{2,\eps,\pm}_{ts}(1,2)|^2 \le C|t-s|^{4\alpha}.$
\end{itemize}

\end{Lemma}

{\bf Remarks.}

\begin{enumerate}
\item Only the {\em increment} $\del G^{2,\eps,\pm}(1,2)$ makes sense: $G^{2,\eps,\pm}(1,2)$ defined as the
integral of $G^{\pm}_t(\xi_1,\xi_2)$ is infra-red divergent (see also remark after Lemma \ref{lem:Kah} below).
\item Even if $\alpha>1/4$, only the {\em sum of the two increment terms} $\del G^{2,\eps,+}_{ts}(1,2)+
\del G^{2,\eps,-}_{ts}(1,2)$ satisfies the above rate of convergence estimate (ii). This is due to the 
spurious singularity on the diagonal $\xi_1=-\xi_2$, not to an ultra-violet divergence when $|\xi_1|,|\xi_2|\to\infty$.
We skip the proof which is
not needed.
\end{enumerate}

{\bf Proof.}  Let us first remark that the  symmetry $(W(1)\leftrightarrow W(2),s\leftrightarrow
t)$ exchanges $G^+_t(\xi_1,\xi_2)-G^+_s(\xi_1,\xi_2)$ with $G^-_t(\xi_2,\xi_1)-G^-_s(\xi_2,\xi_1)$, and
$I_{ts}^+(\xi_1,\xi_2)(\partial)$ with $I_{ts}^-(\xi_2,\xi_1)(\partial)$. Hence it is enough to prove H\"olderianity
for ${\bf B}^{2,\eps,+}(1,2)(\partial)$ and $\del G^{2,\eps,+}(1,2)$.

\medskip 

 We shall use a number of times the following elementary lemma, inspired by arguments of J.-P. Kahane concerning the
regularity of random Fourier series \cite{Kah}:

\begin{Lemma} \label{lem:Kah}

\begin{itemize}

\item[(i)] Let $F(u)=\int_{\R} dW_{\xi} a(\xi)e^{\II u\xi}$, where $|a(\xi)|^2\le C|\xi|^{-1-2\beta}$
for some $0<\beta<1$: then, for every $u_1,u_2\in\R$, 
\BEQ \esper |F(u_1)-F(u_2)|^2 \le C' |u_1-u_2|^{2\beta}.\EEQ

\item[(ii)] Let $\tilde{F}(\eps)=\int_{\R} dW_{\xi} a(\xi)e^{-\eps|\xi|}$ $(\eps>0)$,
 where $|a(\xi)|^2\le C|\xi|^{-1-2\beta}$
for some $0<\beta<1$: then, for every $\eps_1,\eps_2\in\R_+$, 
\BEQ \esper |\tilde{F}(\eps_1)-\tilde{F}(\eps_2)|^2 \le C' |\eps_1-\eps_2|^{2\beta}.\EEQ

\end{itemize}

\end{Lemma}

{\bf Proof.} Bound $|e^{\II u_1\xi}-e^{\II u_2\xi}|$ by $|u_1-u_2| |\xi|$ for $|\xi|\le\frac{1}{|u_1-u_2|}$ and by $2$ otherwise, and similarly for $|e^{-\eps_1|\xi|}-e^{-\eps_2|\xi|}|$.
 Note the variance integral is infra-red convergent near $\xi=0$. 
\hfill \eop

{\bf Remark:} Unless $|a(\xi)|^2$ is $L^1_{loc}$ near $\xi=0$, only the {\em increments} $F(u_1)-F(u_2)$,
$\tilde{F}(\eps_1)-\tilde{F}(\eps_2)$ are well-defined.

\begin{itemize}
\item[(i)] Apply Lemma \ref{lem:Kah} (i)  to $F_s(u)=\int_{\R} dW_{\xi_1}(1) a(\xi_1)e^{\II u\xi_1}$ with
 \BEQ a(\xi_1)=e^{-\eps|\xi_1|} |\xi_1|^{-\half-\alpha} \int_{|\xi_2|\ge |\xi_1|} dW_{\xi_2}(2)
e^{-\eps|\xi_2|} e^{\II s\xi_2} |\xi_2|^{-\half-\alpha}; \EEQ
 since $\Var\  a(\xi_1)\le C|\xi_1|^{-1-4\alpha}$, one gets the uniform H\"olderianity
estimates (i) $\esper ({\bf B}_{ts}^{2,\eps,+}(1,2)(\partial))^2 \le C|t-s|^{4\alpha}$
  for ${\bf B}_{ts}^{2,\eps,+}(1,2)(\partial).$ As for the rate of convergence (ii) (see above 
Lemma \ref{lem:B2}), one rewrites
$e^{-\eps(|\xi_1|+|\xi_2|)}-e^{-\eta(|\xi_1|+|\xi_2|)}$ as $(e^{-\eps|\xi_1|}-e^{-\eta|\xi_1|})
e^{-\eps|\xi_2|}+e^{-\eta|\xi_1|}(e^{-\eps|\xi_2|}-e^{-\eta|\xi_2|})$.  The first term may be bounded
as $F_s$ by applying Lemma \ref{lem:Kah} (ii). For the second term, Lemma \ref{lem:Kah} (ii) should
be applied to $\tilde{F}_s(\eps)=\int_{|\xi_2|\ge |\xi_1|} dW_{\xi_2}(2) e^{-\eps|\xi_2|} e^{\II s\xi_2}
|\xi_2|^{-\half-\alpha}$, which yields an exponent
$2\alpha$ instead of $4\alpha$.

\item[(ii)] Apply Lemma \ref{lem:Kah} (i)  to $F(u)=\int_{\R} dW_{\xi} a(\xi) e^{\II u\xi}$ with
(setting $\xi=\xi_1+\xi_2$)
\BEQ a(\xi)=\frac{1}{\II\xi} \int_{|\xi_2|\ge |\xi-\xi_2|} dW_{\xi_2} |\xi-\xi_2|^{\half-\alpha}
\frac{|\xi_2|^{\half-\alpha}}{\II\xi_2} e^{-\eps(|\xi_2|+|\xi-\xi_2|)}.\EEQ
Setting $\eps$ directly to $0$ yields (assuming for instance $\xi>0$) 
\BEQ \Var\  a(\xi)\le \frac{1}{\xi^2} \int_{\xi/2}^{+\infty} |\xi_2|^{-1-2\alpha} |\xi-\xi_2|^{1-2\alpha} d\xi_2,\EEQ
which converges if and only if $\alpha>1/4$, in which case 
\BEQ \Var\  a(\xi)\le |\xi|^{-1-4\alpha} \int_{\half}^{+\infty} |u|^{-1-2\alpha} |1-u|^{1-2\alpha}
du=C|\xi|^{-1-4\alpha}. \EEQ

\end{itemize}
\hfill \eop

All together one has proved (up to the proof for the rate of convergence for $\del G^{2,\eps,\pm}$) that 
\BEQ {\bf B}_{ts}^{2,\eps}=({\bf B}_{ts}^{2,\eps,+}(1,2)(\partial)+
\del G_{ts}^{2,\eps,+}(1,2))+({\bf B}_{ts}^{2,\eps,-}(1,2)(\partial)+
\del G_{ts}^{2,\eps,-}(1,2)) \EEQ
is $2\alpha^-$-H\"older uniformly in $\eps$, {\it provided} $\alpha>1/4$.

\bigskip

So what should one do when $\alpha\le 1/4$ ?

\medskip

\begin{Definition}[cut Fourier domain]

Let, for some constant $C_{reg}\in(0,1)$, 
 \BEQ \R^2_{reg}:=\{(\xi_1,\xi_2)\ |\ |\xi_1|\le |\xi_2|, |\xi_1+\xi_2|> C_{reg} |\xi_2|\}. \EEQ

\end{Definition}

The condition $|\xi_1+\xi_2|>C_{reg}|\xi_2|$ excludes a conical region along  the singular line
$\xi_1=-\xi_2$.

\begin{Definition}[regularized L\'evy area $ {\cal R} {\bf B}_{ts}^{2,\eps}$ ]

Let
\BEQ {\cal R} {\bf B}_{ts}^{2,\eps}:=({\bf B}_{ts}^{2,\eps,+}(1,2)(\partial)+
\del {\cal R} G_{ts}^{2,\eps,+}(1,2))+({\bf B}_{ts}^{2,\eps,-}(1,2)(\partial)+
\del {\cal R} G_{ts}^{2,\eps,-}(1,2)) \EEQ
where the following {\bf regularized increment term} has been introduced,
\BEQ {\cal R} G_t^{2,\eps,+}(1,2):=c_{\alpha}^2\int\int_{(\xi_1,\xi_2)\in \R^2_{reg}}
  e^{-\eps(|\xi_1|+|\xi_2|)} |\xi_1 \xi_2|^{\half-\alpha} G_t^+(\xi_1,\xi_2)\EEQ
and
\BEQ {\cal R} G_t^{2,\eps,-}(1,2):=c_{\alpha}^2\int\int_{(\xi_2,\xi_1)\in \R^2_{reg}}
  e^{-\eps(|\xi_1|+|\xi_2|)} |\xi_1 \xi_2|^{\half-\alpha} G_t^-(\xi_2,\xi_1).\EEQ
\end{Definition}

The regularized L\'evy area satisfies the multiplicative property (ii) of the Introduction because
the corresponding counterterm 
\BEQ \del{\bf Z}^{2,\eps}(1,2):=\left( \del {\cal R} G^{2,\eps,+}(1,2)-\del G^{2,\eps,+}(1,2) \right)+
\left(  \del {\cal R} G^{2,\eps,-}(1,2)-\del G^{2,\eps,-}(1,2) \right) \EEQ
(given by an integral on the conical Fourier domain $\R^2\setminus\R^2_{reg}$ along the diagonal $\xi_1=-\xi_2$)
is an increment (see Introduction). It satisfies the geometric property (iii) because $\del{\bf Z}^{2,\eps}(1,2)$
is antisymmetric in $1 \leftrightarrow 2$,
 which follows in turn from the symmetry exhibited at the beginning of the proof
of Lemma \ref{lem:B2}.

\medskip

\begin{Theorem}
 For every $\alpha\in(0,1/2)$, ${\cal R} {\bf B}_{ts}^{2,\eps}(1,2)$ is $2\alpha^-$-H\"older uniformly in $\eps$.
\end{Theorem}

{\bf Proof.} Let us first prove the H\"older estimates (i) $\esper |{\cal R} {\bf B}_{ts}^{2,\eps}(1,2)|^2
\le C|t-s|^{4\alpha}$. Similarly to the proof of point (ii) in Lemma \ref{lem:B2}, we apply Lemma \ref{lem:Kah} (i)
to  $F_{reg}(u)=\int_{\R} dW_{\xi} a_{reg}(\xi) e^{\II u\xi}$ with
(setting $\xi=\xi_1+\xi_2$)
\BEQ a_{reg}(\xi)=\frac{1}{\II\xi} \int_{D} dW_{\xi_2} |\xi-\xi_2|^{\half-\alpha}
\frac{|\xi_2|^{\half-\alpha}}{\II\xi_2} e^{-\eps(|\xi_2|+|\xi-\xi_2|)},\EEQ
where $D=\{\xi_2\in\R\ |\ |\xi_2|\ge |\xi-\xi_2|$, $|\xi|>C_{reg}|\xi_2|\}$. In particular,
$|\xi-\xi_2|<\frac{1}{C_{reg}}|\xi|$, which implies $|\xi_2|\le C|\xi|$, so (assuming
for instance $\xi>0$)
\BEA \Var\  a_{reg}(\xi) &\le&  \frac{1}{\xi^2} \int_{\xi/2}^{C\xi} |\xi_2|^{-1-2\alpha} |\xi-\xi_2|^{1-2\alpha} d\xi_2 \nonumber\\
&=& |\xi|^{-1-4\alpha} \int_{1/2}^C |u|^{1-2\alpha} |1-u|^{-1-2\alpha}\ du =C'|\xi|^{-1-4\alpha}<\infty.
\nonumber\\
\EEA

\medskip

For the rate of convergence, we rewrite $e^{-\eps(|\xi_1|+|\xi_2|)}-e^{-\eta(|\xi_1|+|\xi_2|)}$ as the sum of
two terms as
in the proof of Lemma \ref{lem:B2} (i).
 The {\em second} term may be studied  using Lemma \ref{lem:Kah} (ii)  by considering the coefficient
of $dW_{\xi_2}(2)$; since $\frac{|\xi|}{2}\le|\xi_2|<\frac{|\xi|}{C_{reg}}$, one gets the same estimates as for
the proof of the H\"older estimates (i), 
and hence a rate of convergence with exponent $4\alpha$. The {\em first} term involves the function $\tilde{F}(\eps)=\int dW_{\xi_1}
\tilde{a}_{\eta}(\xi_1) e^{\II u\xi_1} e^{-\eps|\xi_1|}$, with
\BEQ \tilde{a}_{\eta}(\xi_1)=|\xi_1|^{\half-\alpha} \int_{|\xi_2|\ge |\xi_1|} dW_{\xi_2} 
\frac{|\xi_2|^{\half-\alpha}}{[\II(\xi_1+\xi_2)][\II\xi_2]} e^{-\eta|\xi_2|}.\EEQ The estimates
$\esper |\tilde{a}_{\eta}(\xi_1)|^2 \le C|\xi|^{-1-4\alpha}$
 is easy to get using the fact that $|\xi_1+\xi_2|>C_{reg}|\xi_2|\ge C_{reg}|\xi_1|$,
which gives once again the same rate of convergence.

 \hfill \eop

{\bf Remark.} Since the boundary terms ${\bf B}^{2,\eps,\pm}_{ts}(1,2)(\partial)$ are uniformly $2\alpha^-$-H\"older
in $\eps$, one could simply have set ${\cal R}G^{2,\eps,\pm}_{ts}(1,2)(\partial)\equiv 0$ (which looks somewhat
drastic). Our point of view is to try and regularize as little as possible.


\section{Iterated integrals of higher order: a sketchy overview}


Combinatorics of 
Fourier normal ordering become non trivial starting from iterated integrals of third order;
we shall concentrate on this case in this section, although some notions will be presented for
the case of general iterated integrals. Detailed proofs should be found in \cite{Unt09}.

Recall the two decompositions of $I_{ts}(\xi_1,\xi_2)$ into $\del G_{ts}^+(\xi_1,\xi_2)
+I_{ts}^+(\xi_1,\xi_2)(\partial)$ and  $\del G_{ts}^-(\xi_2,\xi_1)
+I_{ts}^-(\xi_2,\xi_1)(\partial)$. As already noted, 
the symmetry $\xi_1\leftrightarrow \xi_2$, $s\leftrightarrow
t$ maps $G^+$ to $G^-$ and $I^+(\partial)$ to $I^-(\partial)$, hence terms with indices $\pm$ may
be treated on an equal footing. This is no more the case for the integrals involved in ${\bf B}_{ts}^{3,\eps}$. Namely, Fubini's theorem implies for instance (considering three among the
six permutations of $\{1,2,3\}$, including the trivial one)

\BEQ \int_s^t dB_{u_1}^{\eps}(i_1)\int_s^{u_1}dB_{u_2}^{\eps}(i_2)\int_s^{u_2} dB_{u_3}^{\eps}(i_3)
 \label{eq:3-1} \EEQ
\BEQ  =\int_s^t dB_{u_2}^{\eps}(i_2)\int_{u_2}^{t}dB_{u_1}^{\eps}(i_1)\int_s^{u_2} dB_{u_3}^{\eps}(i_3)
 \label{eq:3-2} \EEQ
\BEQ  \qquad\qquad =\int_s^t dB_{u_2}^{\eps}(i_2)\int_s^{u_2}dB_{u_3}^{\eps}(i_3)\int_{u_2}^{t} dB_{u_1}^{\eps}(i_1)  \label{eq:3-3}
\EEQ

Each of these may be represented as a finite sum of {\it tree iterated integrals} as we shall
presently see. If $\T$ is a decorated rooted tree, i.e. a tree  with a distinguished vertex $v_1$
 called {\it root}, such that each vertex $v\in V(\T)=\{$vertices of $\T\}$
wears a label $\ell(v)\in\{1,\ldots,d\}$, and $\Gamma=(\Gamma(1),\ldots,\Gamma(d))$ is a smooth
$d$-dimensional path, then the integral of $\Gamma$ along $\T$ is (denoting by $v^-$, $v\in V(\T)\setminus\{
v_0\}$ the unique ancestor of $v$, i.e. the unique vertex just below $v$)
\BEQ [I_{\T}(\Gamma)]_{ts}:=\int_s^t d\Gamma_{x_{v_1}}(\ell(v_1))\int_s^{x_{v_2^-}} d\Gamma_{x_{v_2}}(\ell(v_2))
\ldots \int_s^{x_{v^-_{|V(\T)|}}} d\Gamma_{x_{v_{|V(\T)|}}}(\ell(v_{|V(\T)|}))
\EEQ
where $(v_1,\ldots,v_{|V(\T)|})$ is any ordering of $V(\T)$ compatible with the tree partial ordering, i.e. such that $v^-<v$ for every $v\in V(\T)\setminus\{v_1\}$
 (in other words, such that the indices of the vertices decrease while going down the branches towards the root). The definition extends easily to {\it
forests}, i.e. to (finite) disjoint unions of trees (which are seen as a commutative product of the
trees),
 by multiplying the tree iterated integrals
corresponding to each connected component of the forest, and then (by taking linear combinations)
to the {\bf algebra over $\R$ generated by decorated rooted trees}, $\cal T$, which is actually a {\bf Hopf algebra} \cite{ConKre98,ConKre00}.

\begin{figure}[h]
  \centering
   \includegraphics[scale=0.60]{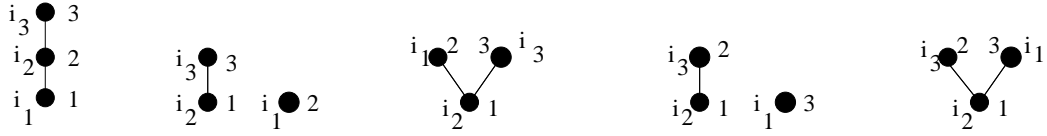}
   \caption{\small{Example of iterated integrals. From left to right: $\T_1, \T_{2,1},
\T_{2,2},\T_{3,1},\T_{3,2}$.}}
 \label{Fig1}
\end{figure}

Now replace in eq. (\ref{eq:3-1},\ref{eq:3-2},\ref{eq:3-3}) $\int_u^t$ by $\int_s^t -\int_s^u$ and
$\int_u^{u'}$ by $\int_s^{u'}-\int_s^u$. Then eq.  (\ref{eq:3-1},\ref{eq:3-2},\ref{eq:3-3}) may
be represented resp. as $I_{\T_1}(B^{\eps})$, $I_{\T_{2,1}}(B^{\eps})-I_{\T_{2,2}}(B^{\eps})$,
$I_{\T_{3,1}}(B^{\eps})-I_{\T_{3,2}}(B^{\eps})$ (see Fig. \ref{Fig1}).

More generally, to each permutation $\sigma$ of $\{1,2,3\}$ corresponds a rewriting of
 ${\bf B}_{ts}^{3,\eps}(i_1,i_2,i_3)$ as some finite sum,
\BEQ {\bf B}_{ts}^{3,\eps}(i_1,i_2,i_3)=\sum_j g(\sigma,j) I_{\T_j^{\sigma}}(B^{\eps}), \label{eq:g} \EEQ
where $g(\sigma,j)=\pm 1$ is a sign.
Note that each $I_{\T_j^{\sigma}}(B^{\eps})$ equals more explicitly (after permutation of the
names of the $\xi$-variables, $\xi_i\to\xi_{\sigma(i)}$)
\BEA && [I_{\T_j^{\sigma}}(B^{\eps})]_{ts}=c_{\alpha}^3 \int_{\R} dW_{\xi_1}(i_{\sigma(1)})
\int_{\R} dW_{\xi_2}(i_{\sigma(2)}) \int_{\R} dW_{\xi_3}(i_{\sigma(3)}) |\xi_1\xi_2\xi_3|^{\half-
\alpha} \ .\nonumber\\
&& \qquad e^{-\eps(|\xi_1|+|\xi_2|+|\xi_3|)} \int_s^{x_{v_1^-}} e^{\II x_{v_1}\xi_1} dx_{v_1} \int_s^{x_{v_2^-}} e^{\II x_{v_2}\xi_2}
dx_{v_2} \int_s^{x_{v_3^-}}
e^{\II x_{v_3}\xi_3} dx_{v_3}, \label{eq:ITjsigma} \nonumber\\  \EEA
where $(v_1,v_2,v_3)$ is the natural ordering of $V(\T_j^{\sigma})$ (compatible with the tree ordering) given
by the position of the corresponding variable of integration inside the iterated integral {\it after}
applying  Fubini's theorem, see eq. (\ref{eq:3-1},\ref{eq:3-2},\ref{eq:3-3}) or Figure \ref{Fig1}, and
$x_{v_i^-}=t$ if $v_i$ is a root. Note also that the labels on the trees are simply $\ell(v_j)=i_{\sigma(j)}$.

\medskip

The idea of Fourier normal ordering consists in rewriting first ${\bf B}_{ts}^{3,\eps}(i_1,i_2,i_3)$
as (letting $\Sigma_3$ be the group of permutations of $\{1,2,3\}$)
\BEA &&
c_{\alpha}^3 \sum_{\sigma\in\Sigma_3} \sum_j g(\sigma,j) \int\int\int_{|\xi_1|\le|\xi_2|\le|\xi_3|} dW_{\xi_{v_1}}(i_{\sigma(1)}) dW_{\xi_{v_2}}(i_{\sigma(2)}) dW_{\xi_{v_3}}(i_{\sigma(3)}) \ . \nonumber\\
&&  \qquad  .\  e^{-\eps(|\xi_1|+|\xi_2|+|\xi_3|)}  \int_s^{x_{v_1^-}} e^{\II x_{v_1}\xi_1} dx_{v_1} \int_s^{x_{v_2^-}} e^{\II x_{v_2}\xi_2}
dx_{v_2} \int_s^{v_3^-}
e^{\II x_{v_3}\xi_3} dx_{v_3}, \nonumber\\ \EEA
where $\{v_1,v_2,v_3\}$ are the vertices of $\T^{\sigma}_j$, 
so that innermost integrals bear highest Fourier indices.

\medskip

The next task is to get rid of divergences by adding some counterterms (or in other words, by
 discarding
the contribution to the  integral of adequate Fourier subdomains). If one wants the multiplicative
property to remain true, this must be done in compatibility with lower-order counterterms (i.e. with
the definition of the regularized L\'evy area at this stage). Hence one should have:
\BEQ {\bf B}_{ts}^{3,\eps}(i_1,i_2,i_3)-{\bf B}_{tu}^{3,\eps}(i_1,i_2,i_3)-{\bf B}_{us}^{3,\eps}(i_1,i_2,i_3)={\bf B}_{tu}^{1,\eps}(i_1) {\bf B}^{2,\eps}_{us}(i_2,i_3)+{\bf B}_{tu}^{2,\eps}(i_1,i_2)
{\bf B}_{us}^{1,\eps}(i_3). \EEQ

This identity has been generalized to tree integrals (see \cite{Gu2}). Namely, letting $\T$ be
a tree and $\vec{v}$ range over all {\em admissible cuts} of $\T$ (we shall use the notation:
$\vec{v}\models V(\T)$), i.e. over all
non-empty subsets $\{v_1,\ldots,v_J\}\in V(\T)\setminus\{0\}$ ($0$=root of $\T$)
 such that no pair $\{v_i,v_j\}\subset\vec{v}$ is connected by going down or up the tree, then

\begin{figure}[h]
  \centering
   \includegraphics[scale=0.70]{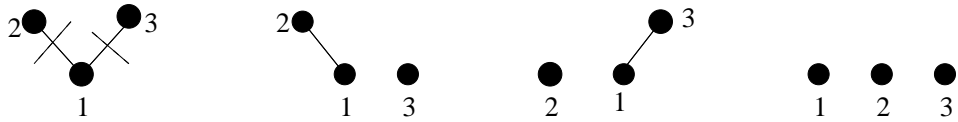}
   \caption{\small{Admissible cuts of $\T_{2,2}$: $\vec{v}=\{3\},\{2\}$ or $\{2,3\}$ is also in this particular
case
the set of vertices of $R_{\vec{v}}\T$.}}
 \label{Fig2}
\end{figure}

\BEQ [\del I_{\T}(B^{\eps})]_{tus}=\sum_{\vec{v}\models V(\T)} [I_{L_{\vec{v}}\T}(B^{\eps})]_{tu}
[I_{R_{\vec{v}} \T}(B^{\eps})]_{us} \label{eq:tree-x} \EEQ

where $R_{\vec{v}}\T\subset\T$ is the forest obtained as the union of branches lying above
$v_1,\ldots,v_J$ (including the vertices $v_1,\ldots,v_J$), and $L_{\vec{v}}\T\subset\T$ is the subtree obtained after removing these branches
(see example on Fig. \ref{Fig2}). This identity called {\bf tree multiplicative property} may be rephrased \cite{Unt09bis} by saying that $[\del I(B^{\eps})]_{tus}$ (viewed as a linear form on
the Hopf algebra $\cal T$) is the {\em convolution} of $[I(B^{\eps})]_{tu}$ and $[I(B^{\eps})]_{us}$.

\medskip

Now (see footnote  after Lemma \ref{lem:B2}), one wants to set formally $s=\pm\II\infty$ to generalize the
increment/boundary decomposition of the previous section. This results in the following definition
of {\it skeleton integrals}:

\begin{Definition}[skeleton integrals]

Let 
\BEA && [\Sk I_{\T}(B^{\eps})]_{t}=c_{\alpha}^3 \int_{\R} dW_{\xi_1}(i_{\sigma(1)})
\int_{\R} dW_{\xi_2}(i_{\sigma(2)}) \int_{\R} dW_{\xi_3}(i_{\sigma(3)}) |\xi_1\xi_2\xi_3|^{\half-
\alpha} \ .\nonumber\\
&& \qquad e^{-\eps(|\xi_1|+|\xi_2|+|\xi_3|)} \int^{x_{v_1^-}} e^{\II x_{v_1}\xi_1} dx_{v_1} \int^{x_{v_2^-}} e^{\II x_{v_2}\xi_2} dx_{v_2}
 \int^{x_{v_3^-}}
e^{\II x_{v_3}\xi_3} dx_{v_3} \nonumber\\ \EEA
(see eq. (\ref{eq:ITjsigma})).

\end{Definition}

 Formally $[\Sk I_{\T}(B^{\eps})]_t=[I_{\T}(B^{\eps})]_{t,\pm\II\infty}.$ As was the case for $G^{2,\eps,\pm}(1,2)$, only the {\it increment} $[\del\Sk I_{\T}(B^{\eps})]_{ts}=[\Sk I_{\T}(B^{\eps})]_t -
 [\Sk I_{\T}(B^{\eps})]_s$ makes sense.
The tree multiplicative property yields then the {\bf tree skeleton decomposition}
\BEQ [I_{\T}(B^{\eps})]_{tu}=[\Sk I_{\T}(B^{\eps})]_t-[\Sk I_{\T}(B^{\eps})]_u-
\sum_{\vec{v}\models V(\T)} [I_{L_{\vec{v}}\T}(B^{\eps})]_{tu} [\Sk I_{R_{\vec{v}}\T}(B^{\eps})]_u.
\\ \label{eq:sk}
\EEQ

This is an inductive formula, which yields $I_{\T}(B^{\eps})$ in terms of integrals or skeleton integrals of lower order. Once again, it may be interpreted as a convolution of characters of the Hopf
algebra $\cal T$, as the convolution of $[\Sk I(B^{\eps})]_t$ with the inverse of  $[\Sk I(B^{\eps})]_u$
(defined via the {\em antipode}) to be precise, see \cite{Unt09bis}. If e.g. $\T$ has 3 vertices, then each tree component of $R_{\vec{v}}\T$ has
at most 2 vertices. A skeleton integral for a tree reduced to one vertex is simply
\BEQ \int_{\R} dW_{\xi} e^{-\eps|\xi|} |\xi|^{\half-\alpha} \int^u e^{\II u\xi} du=
\int_{\R} dW_{\xi} e^{-\eps|\xi|} e^{\II u\xi} \frac{|\xi|^{\half-\alpha}}{\II\xi}\EEQ
(whose increments are those of  $B^{\eps}$ itself) and requires no regularization,
whereas for a tree with two vertices, one gets
\BEA &&\int_{\R} dW_{\xi_1}(1)\int_{\R} dW_{\xi_2}(2) e^{-\eps(|\xi_1|+|\xi_2|)}
|\xi_1\xi_2|^{\half-\alpha} \int^u e^{\II x_{v_1}\xi_1} dx_{v_1}\int^{x_{v_1}} e^{\II x_{v_2}\xi_2} dx_{v_2} \nonumber\\
&&= \int_{\R} dW_{\xi_1}(1)\int_{\R} dW_{\xi_2}(2) e^{-\eps(|\xi_1|+|\xi_2|)}
|\xi_1\xi_2|^{\half-\alpha} G_u^+(\xi_1,\xi_2)\EEA
(which is divergent in the limit $\eps\to 0$ when $\alpha<1/4$, as proved in Lemma \ref{lem:B2}).

The reader may easily check that formula (\ref{eq:sk}) gives precisely the increment/boundary
decomposition of the previous section when $\T$ has two vertices.

\medskip

In order to get convergent quantities in compatibility with the regularization of second-order
integrals, one should (i) regularize the new skeleton integrals of third order
 $[\Sk I_{\T}(B^{\eps})]_t$; (ii) replace $[I_{L_{\vec{v}}\T}(B^{\eps})]_{tu}$ and 
$[\Sk I_{R_{\vec{v}}\T}(B^{\eps})]_u$ in the right-hand
side of (\ref{eq:sk}) by the corresponding regularized quantity of order (1 or) 2.

\medskip

In other words, {\bf regularization must be performed on each  skeleton integral of order $\ge 2$}.

\begin{Definition}[regularized skeleton integrals of order 3]

Let, for $\T=\T_1,\T_{2,1},\T_{2,2},\T_{3,1}$ or $\T_{3,2}$ corresponding to one of the three
above permutations $\sigma\subset\Sigma_3$,
\BEA && [{\cal R}\Sk I_{\T}(B^{\eps})]_u=c_{\alpha}^3\int\int\int_{(\xi_1,\xi_2,\xi_3)\in\R^{\T}_{reg}}
dW_{\xi_1}(i_{\sigma(1)})dW_{\xi_2}(i_{\sigma(2)})dW_{\xi_3}(i_{\sigma(3)}) \nonumber\\
&& \qquad  e^{-\eps(|\xi_1|+
|\xi_2|+|\xi_3|)} |\xi_1 \xi_2\xi_3|^{\half-\alpha} \int^{x_{v_1^-}} e^{\II x_{v_1}\xi_1} dx_{v_1} \int^{x_{v_2^-}} e^{\II x_{v_2}\xi_2}  dx_{v_2} \int^{x_{v_3^-}} e^{\II x_{v_3}\xi_3} dx_{v_3},\nonumber\\ \EEA
where $\R^{\T}_{reg}\subset\R^3$ is the subdomain of integration defined by:
\begin{itemize}
\item[(i)] $\R^{\T_1}_{reg}=\{|\xi_1|\le |\xi_2|\le |\xi_3|;\ \ |\xi_2+\xi_3|,|\xi_1+\xi_2+\xi_3|>
C'_{reg} |\xi_3|\};$
\item[(ii)] $\R^{\T_{2,1}}_{reg}=\{ |\xi_1|\le |\xi_2|\le |\xi_3|;\ \  |\xi_1+\xi_3|>C'_{reg}
|\xi_3|\};$
\item[(iii)] $\R^{\T_{2,2}}_{reg}=\{ |\xi_1|\le |\xi_2|\le |\xi_3|;\ \  |\xi_1+\xi_2+\xi_3|>C'_{reg}
|\xi_3|\};$
\item[(iv)] $\R^{\T_{3,1}}_{reg}=\{ |\xi_1|\le |\xi_2|\le |\xi_3|;\ \  |\xi_1+\xi_2|>C'_{reg}
|\xi_2|\};$
\item[(v)] $\R^{\T_{3,2}}_{reg}=\{ |\xi_1|\le |\xi_2|\le |\xi_3|;\ \  |\xi_1+\xi_2+\xi_3|>C'_{reg}
|\xi_3|\}$
\end{itemize}

for some constant $C'_{reg}\in(0,1)$.

\end{Definition}

For the general definition of the cut Fourier domains of integrations, we refer the reader to
\cite{Unt09}. Uniform H\"olderianity with respect to $\eps$ has been defined before Lemma \ref{lem:B2}.

\begin{Theorem}

Assume $\alpha<1/3$. Then 
the above regularized skeleton integrals are $3\alpha^-$-H\"older uniformly in $\eps$.
\end{Theorem}

{\bf Proof} (sketch). Let us just sketch the proof e.g. for $[{\cal R}\Sk I_{\T_1}(B^{\eps})]_u$, assuming
$i_1\not=i_2\not=i_3$, say, $i_1=1$, $i_2=2$, $i_3=3$.
The triple integral $\int^{x_{v_1^-}}  e^{\II x_{v_1}\xi_1} dx_{v_1} \int^{x_{v_2^-}} e^{\II x_{v_2}\xi_2} 
dx_{v_2} \int^{x_{v_3^-}} e^{\II x_{v_3}\xi_3} dx_{v_3}$ writes (by straightforward computation) $\frac{e^{\II t(\xi_1+\xi_2+\xi_3)}}{
[\II \xi_3][\II(\xi_2+\xi_3)][\II(\xi_1+\xi_2+\xi_3)]}$. The idea is that the restriction of the domain of
integration ensures that denominators are not {\it too} small. Apply Lemma \ref{lem:Kah} (i) to
$F(u)=\int_{\R} dW_{\xi} a(\xi) e^{\II u\xi}$ with (setting $\xi=\xi_1+\xi_2+\xi_3$)
\BEQ a(\xi)=\frac{1}{\II\xi} \int_{D_{\xi}} dW_{\xi_2}(2)dW_{\xi_3}(3) |\xi-\xi_2-\xi_3|^{\half-\alpha} \frac{|\xi_2|^{\half-\alpha}}{\II(\xi_2+\xi_3)} \frac{|\xi_3|^{\half-\alpha}}{\II\xi_3} \EEQ
on a domain $D_{\xi}$ on which $\frac{|\xi|}{3}\le |\xi_3|<\frac{|\xi|}{C'_{reg}}$, $|\xi-\xi_2-\xi_3|^{\half-\alpha}\le
|\xi_3|^{\half-\alpha}$ and $\frac{|\xi_2|^{\half-\alpha}}{|\xi_2+\xi_3|}\le \frac{1}{C'_{reg}}
|\xi_3|^{-\half-\alpha}$. Considering $\Var \ a(\xi)\le \frac{C}{\xi^2} \int\int_{D_{\xi}}
d\xi_2 d\xi_3 |\xi_3|^{-1-6\alpha}$, the integral over $|\xi_2|\le |\xi_3|$ contributes $O(|\xi_3|)$,
while the integral over $\frac{|\xi|}{3}\le |\xi_3|<\frac{|\xi|}{C'_{reg}}$ leads to $\frac{C}{\xi^2}O(|\xi|^{1-6\alpha})
=O(|\xi|^{-1-6\alpha}).$ Hence $\esper \left( [\del {\cal R}\Sk I_{\T_1} (B^{\eps})]_{ts} \right)^2\le
C|t-s|^{6\alpha}.$  By applying Lemma \ref{lem:Kah} (ii), a similar convergence rate may be
proved. \hfill \eop

\medskip

We may now finally define regularized iterated integrals of order 3 according to the above scheme.

\begin{Definition}[regularized integrals of order 3]

Let 
\BEQ [{\cal R}I_{\T}({\bf B}^{\eps})]_{tu}(i_1,i_2,i_3)=\sum_{\sigma\in\Sigma_3}
\sum_j g(\sigma,j) [{\cal R}I_{\T^{\sigma}_j}(B^{\eps})]_{tu} \EEQ where $g(\sigma,j)$ is
as in eq. (\ref{eq:g}), and by definition
\BEQ [{\cal R}I_{\T}(B^{\eps})]_{tu}=[{\cal R}\Sk I_{\T}(B^{\eps})]_t-[{\cal R}\Sk I_{\T}(B^{\eps})]_u-
\sum_{\vec{v}\models V(\T)} [{\cal R}I_{L_{\vec{v}}\T}(B^{\eps})]_{tu} [{\cal R}\Sk I_{R_{\vec{v}}\T}(B^{\eps})]_u.
\\ \label{eq:Rsk}
\EEQ

\end{Definition}

The definition mimicks the previous unregularized tree multiplicative property for skeleton
integrals, eq. (\ref{eq:sk}).
All terms in it are convergent when $\eps\to 0$. There remains only to prove that the multiplicative/Chen (ii)
and geometric/shuffle property (iii) of the Introduction  are preserved by the regularization. In principle (ii)
 should be true because
eq. (\ref{eq:Rsk}) is a multiplicative property in itself; it must only be shown that this multiplicative
property implies the original multiplicative property (ii) of the Introduction after summing over all
permutations $\sigma$, which is best done in a Hopf algebra language.
 The geometric/shuffle property is proved by showing that the regularized  integration operator
$\sum_{\sigma} \sum_j g(\sigma,j) {\cal R}\Sk I_{\T^{\sigma}_j}$ is a character of another Hopf algebra called {\em
shuffle algebra}. Proofs should be found in  \cite{Unt09bis}.


\end{document}